\documentclass[10pt]{article}
\usepackage{indentfirst, latexsym, amsmath,bm,amssymb}
\usepackage{tikz}
\usepackage{color,soul}
\usepackage{amsmath}
\usepackage{amssymb}
\usepackage{mathrsfs}
\usepackage{mathtools}
\usepackage{stmaryrd}
\usepackage{bbold}

\allowdisplaybreaks[4]
\usepackage{amsthm}
\newtheorem{thm}{Theorem}[section]
\newtheorem{lem}[thm]{Lemma}
\newtheorem{coro}[thm]{Corollary}

\newtheorem{defi}[thm]{Definition}
\newtheorem{rem}{Remark}[]
\newenvironment {Proof}{\noindent {\bf Proof.}}{\hfill\ensuremath{\square}}
\newenvironment {Proof of 1}{\noindent {\bf Proof of Theorem \ref{thm1}.}}{\hfill\ensuremath{\square}}
\newenvironment {Proof of 2}{\noindent {\bf Proof of Theorem \ref{thm2}.}}{\hfill\ensuremath{\square}}

\begin{document}

\title{ On weighted spectral radius of unraveled balls and normalized Laplacian eigenvalues \thanks{This work is partly supported by the National Natural Science Foundation of China (Nos. 11971311,12161141003, 12026230) and the
Montenegrin-Chinese Science and Technology Cooperation Project (No.3-12).
}}

\author{ Yuzhenni Wang, Xiao-Dong Zhang\footnote{Corresponding author. E-mail: xiaodong@sjtu.edu.cn (X.-D. Zhang)}\\
School of Mathematical Sciences, MOE-LSC, SHL-MAC\\
Shanghai Jiao Tong University, Shanghai 200240, P. R. China\\
Email: wangyuzhenni@sjtu.edu.cn, xiaodong@sjtu.edu.cn.\\
}

\date{}

\maketitle
\begin{abstract}
  For a graph $G$, the unraveled ball of radius $r$ centered at a vertex $v$ is the ball of radius $r$ centered at $v$ in the universal cover of $G$.   We obtain a lower bound on the weighted spectral radius of unraveled balls of fixed radius in a graph with positive weights on edges, which is used to present an upper bound on the 
   $s$-th (where $s\ge 2$) smallest normalized Laplacian eigenvalue of irregular graphs under minor assumptions. Moreover, when $s=2$, the result may be regarded as an Alon--Boppana type bound for a class of irregular graphs.

\end{abstract}

\section{Introduction}
For a simple graph $G$ of order $n$, denote the eigenvalues of its adjacency matrix by $\lambda_1(G)\ge\lambda_2(G)\ge \dots\ge\lambda_n(G)$. The degree of a vertex $u$ in $G$, denoted by $d_G(u)$, is defined as the number of edges incident to $u$, and we write $d(u)$ when $G$ is clear.

The well-known result of Alon--Boppana bound may be stated as follows:
\begin{thm}\cite{alon1986, Nilli1991}
For any $d$-regular graph $G$ containing two edges with distance at least $ 2k+2$,
\begin{equation}\label{eq-Alon}
 \lambda_2(G)\ge2\sqrt{d-1}\left(1-\frac{1}{k+1}\right)+\frac{1}{k+1}.
\end{equation}
\end{thm}
It is natural to generalize the Alon--Boppana bound to graphs that are not necessarily regular. One plausible extention is to verify whether the inequality $\liminf_{i\rightarrow\infty}\lambda_2(G_i)\ge 2\sqrt{d-1}$ still holds for any sequence of graphs $\{G_i\}$ with average degree at least $d$ and growing diameter. However, Hoory \cite{Hoory2005} constructed a counterexample to disprove the problem. Furthermore, he proposed a new idea of $r$-robust and  provided an Alon--Boppana type bound for a class of irregular graphs as follows.
For a graph $G$, the \emph{ball} of radius $r$ centered at $v$, denoted by $G(v,r)$, is the induced subgraph of $G$ on the vertices with distance at most $ r$ from $v$.
 A graph $G$ \emph{has an $r$-robust average degree $\ge d$} if for the induced subgraph obtained by deleting any ball of radius $r$, its average degree is at least $d$.
\begin{thm}\cite{Hoory2005}\label{thm-Hoory}
Given a real number $d\ge 2$ and a natural number $r\ge 2$, if a graph $G$ has an $r$-robust average degree $\ge d$, then
\begin{equation}\label{eq-Hoory}
 \max\{\lambda_2(G),|\lambda_{n}(G)|\}\ge 2\sqrt{d-1}\left(1-c\cdot \frac{\log r}{r}\right),
\end{equation}
where $c$ is an absolute constant.
\end{thm}
Recently, Jiang \cite{Jiang2019} presented a method of unraveled balls to improve the bound above.
\begin{thm}\cite{Jiang2019}\label{thm-Jiang}
Given a real number $d\ge 1$ and a natural number $r\ge 1$, if a graph $G$ has an $r$-robust average degree $\ge d$, then
\begin{equation}\label{eq-Jiang}
 \lambda_2(G)\ge 2\sqrt{d-1}\cos\left(\frac{\pi}{r+1}\right).
\end{equation}
\end{thm}
It is benefit to consider normalized Laplacian eigenvalues, since they can reveal many fundamental properties of a graph (see\cite{chung1997}). Especially, the second smallest one is tightly related with expansion and algorithmic properties of a graph (see \cite[Chapter 2]{chung1997}).  The normalized Laplacian matrix $\mathcal {L}(G)$ of a graph $G$ is defined to be
 $I-D^{-\frac 12}AD^{-\frac 12}$, where $D$ is the diagonal degree matrix with diagonal element $D(v,v)=d(v)$ for $v\in V(G)$, and $A$ is the adjacency matrix of $G$. Denote the eigenvalues of $\mathcal {L}(G)$ by $0=\mu_1(G)\le\mu_2(G)\le \dots\le\mu_n(G)\le 2$. In terms of the normalized Laplacian, the well-known Alon--Boppana bound theorem says that for any $d$-regular graph $G$,
\begin{equation}\label{eq-Alon-nor}
 \mu_2(G)\le1-2\frac{\sqrt{d-1}}{d}+o(1),
\end{equation}
as the diameter of $G$ going to infinity.
One may ask whether the assertion that $\limsup_{i\rightarrow \infty}\mu_2(G_i)\le1-2\frac{\sqrt{d-1}}{d}+o(1)$ for any sequence of graphs $G_i$ with average degree at least $ d$ and growing diameter still holds.
Indeed, Young\cite{Young2022} proved that there exists some fixed $\varepsilon\ge 0$ and a sequence of graphs $\{G_i\}$ with the common average degree $d$ and common maximum degree (and hence growing diameter) such that for every $i$,
\begin{equation}
\mu_2(G_i)\ge 1-2\frac{\sqrt{d-1}}{d}+\varepsilon.
\end{equation}
 Furthermore, Young\cite{Young2022} generalized the Alon--Boppana bound on the second smallest normalized Laplacian eigenvalue to graphs that may be irregular by adapting the idea of Hoory\cite{Hoory2005}.
Young indicated that \cite[Theorem 7]{Young2022} can be phrased in the $r$-robust average degree framework of Hoory.
The \emph{second order average degree} of a graph $G$ is defined to be
$$\widetilde{d}_{G}=\frac{\sum_{u\in V(G)}d(u)^2}{\sum_{u\in V(G)}d(u)}.$$ A graph $G$ is \emph{$(r,d,\widetilde{d})$-robust} if for the induced subgraph of $G$ obtained by deleting any ball of radius $r$,
its average degree is at least $d$, and its second order average degree is at most $ \widetilde{d}$.
\begin{thm}\cite{Young2022}\label{thm-Young}
Given real numbers $\widetilde{d}\ge d\ge 2$ and a natural number $r\ge 2$, if a graph $G$ is $(r,d,\widetilde{d})$-robust, then
\begin{equation}\label{eq-Young}
 \mu_2(G)\le 1- \frac{2\sqrt{d-1}}{\widetilde{d}} \left(1-c\cdot\frac {\log r}{r}\right).
\end{equation}
\end{thm}

In addition,  Chung \cite{Chung2016} used a different approach to obtain,  under some technical assumptions on a graph, another analogous upper bound $\mu_2(G)\le 1-\sigma\left(1-\frac ck\right)$, where $\sigma=\frac{2\sum_{u\in V (G)} d(u)\sqrt{d(u)-1}}{\sum_{u\in V(G)}d(u)^2}$, $k$ is the diameter of $G$, and $c$ is a constant.

Let $G$ be a simple graph. The matrix $D^{-\frac 12}AD^{-\frac 12}$ associated with $G$ may be regarded as the  adjacency matrix of the weighted graph $(G,w_0)$ with edge weight $w_0(uv)=(d(u) d(v))^{-\frac 12}$.
Moreover,  if we denote the second largest eigenvalue of $D^{-\frac 12}AD^{-\frac 12}$ by $\lambda_2(G,w_0)$, then the second smallest normalized Laplacian eigenvalue $\mu_2(G)$ is equal to
$1-\lambda_2(G,w_0)$.
Based on the observations above, it is believable that by considering weighted graphs one could provide a tighter upper bound on $\mu_2(G)$. The related results are referred to \cite{angel2015, chung1997, Srivastava2018, Mohar2010}. The motivation of this paper is to solve the problem above by combining the weighting idea of Young\cite{Young2022} and the idea of unraveled ball of Jiang\cite{Jiang2019}. Indeed, we present an upper bound on $\mu_s(G)$ for $s\ge 2$.
Hereafter, denote the set of real numbers and positive real numbers by $\mathcal{R}$ and $\mathcal{R}^+$ respectively.

Recall that the degree of a vertex $u$ in $G$, denoted by $d(u)$, is the number of edges incident to $u$. One of the main result in this paper is as follows:
\begin{thm}\label{thm1}
If a connected (positively) weighted graph $(G,w)$ has minimum degree at least 2, then for any natural number $r$ with $r\ge 1$ and any function $g\colon V(G)\rightarrow \mathcal{R}^+$, there exists a vertex $v$ of $G$ such that the weighted spectral radius of the unraveled ball $\widetilde{G}(v,r)$ of $G$ satisfies
\begin{equation}\label{eq-thm1}
 \lambda_1(\widetilde{G}(v,r),w)
 \ge \frac{2\sum_{v_1\in V(G)}\sqrt{d(v_1)-1}\sum_{v_2\in N(v_1)}~w(v_1v_2)~\sqrt{g(v_1)g(v_2)}}{\sum_{v\in V(G)} g(v)d(v)}\cos\left(\frac{\pi}{r+2}\right),
\end{equation}
where the neighborhood $N(v_1)$ of $v_1$ is $\{v\in V(G)\colon v_1v\in E(G)\}$.
\end{thm}
\begin{rem}
\begin{description}
      \item[(1)]
  The function $g(x)$ in Theorem~\ref{thm1} is just a technical function which may be appropriately chosen to get simpler bounds for different weight graphs. For instance, if $w\equiv1$, then let $g\equiv 1$ in Theorem~\ref{thm1}, and we obtain \cite[Theorem 1]{Jiang2019}.
  As a result, Theorem~\ref{thm1} extends the result of Jiang to weighted graphs.
  \item[(2)] A graph is positively weighted if the graph has a positive weight on every edge.
    \end{description}
\end{rem}
The other main result in this paper is to derive an upper bound on the 
 $s$th smallest normalized Laplacian eigenvalue. A graph $G$ is \emph{$(r,d,\widetilde{d},s)$-robust} if for the induced subgraph of $G$ obtained by sequentially deleting any $s$ ball of radius $r$,
its average degree is at least $d$, and its second order average degree is at most $ \widetilde{d}$. When $s=1$,  $(r,d,\widetilde{d},1)$-robust is just $(r,d,\widetilde{d})$-robust. We prove the following:
\begin{thm}\label{thm2}
Given real numbers, $\widetilde{d}\ge d\ge 2$, and natural numbers, $r\ge 1$ and $s\ge 2$, if a graph $G$ is $(r,d,\widetilde{d},s-1)$-robust, then the
 $s$-th smallest normalized Laplacian eigenvalue satisfies
\begin{equation}\label{eq-m2}
 \mu_s(G)\le 1- \frac{2\sqrt{d-1}}{\widetilde{d}}\cos\left(\frac{\pi}{r+1}\right).
\end{equation}
\end{thm}
\begin{rem}
 Note that $$\cos\left(\frac{\pi}{r+1}\right)
= 1-\frac{\pi^2}{2(r+1)^2}+o\left(\frac{1}{r^3}\right).$$
 It is  easy to see that  \eqref{eq-m2} is a slight improvement of \eqref{eq-Young} for large $r$ and $s=2$.
\end{rem}
The rest of the paper is organized as follows. In Section 2, some related concepts and symbols are introduced. In Section 3, we prove Theorem \ref{thm1}, and include some corollaries. In Section 4, we present a lower bound on the weighted spectral radius of a ball, which is used to prove Theorem \ref{thm2} in Section 5.

\section{Preliminary}
\begin{defi}
A graph $G^{\prime}$ (possibly infinite) is a covering of another graph $G$ via a covering map $\varphi\colon V(G^{\prime})\rightarrow V(G)$ if $\varphi$ is a surjective map and is a local isomorphism: for every vertex $v$ of $G^{\prime}$, the map $\varphi$ induces a bijection from the edges incident to $v$ in $G^{\prime}$ to the edges incident to $\varphi(v)$ in $G$.
\end{defi}
\begin{defi}
The universal cover $\widetilde{G}$ of a connected graph $G$ is a covering of $G$ that is a (possibly infinite) tree.
\end{defi}
The universal cover of a connected graph is unique up to isomorphism\cite{Leighton1982,Bordenave2019}. If $G$ is a finite tree, then the universal cover of $G$ is itself. Otherwise, the universal cover of $G$ is an infinite graph.
 For instance, the universal cover of $d$-regular graph is the infinite $d$-regular tree.

A \emph{non-backtracking walk} of $G$ is defined as a walk $(v_0,v_1,\dots)$ on $G$ satisfying $v_i\neq v_{i+2}$, for every $i$ with $i\ge 0$. Specifically, a walk of length at most 1 is just non-backtracking.
From the view of random walks, the universal cover can be defined in an equivalent way:
\begin{defi}\cite{Leighton1982}
 The universal cover $\widetilde{G}$ of a connected graph $G$ is defined as follows: the vertex set consists of all non-backtracking walks on $G$ starting at a fixed vertex $v_0$, and two vertices are adjacent if and only if one is a simple extension of the other. The covering map $\varphi\colon V(\widetilde{G})\rightarrow V(G)$ is defined by $\varphi((v_0,\dots,v_i))=v_i$ for $(v_0,\dots,v_i)\in V(\widetilde{G})$.
\end{defi}
In fact, the universal cover is independent of the choice of the fixed vertex $v_0$.
Given a graph $G$, the unraveled ball of $G$ is the ball of radius $r$ centered at $v$ in the universal cover $\widetilde{G}$ of $G$, and thus it has an equivalent definition:
\begin{defi}
Given a graph $G$ and a vertex $v$ of $G$, the unraveled ball of $G$, denoted by $\widetilde{G}(v,r)$, is defined as follows: the vertex set contains all non-backtracking walks on $G$ of length at most $r$ starting at $v$, and two vertices are adjacent if and only if one is a simple extension of the other.
\end{defi}
Next we introduce weighted graphs.
\begin{defi}
A weighted graph $(G,w)$ is a graph $G$ along with a weight function on edges, $w\colon E(G)\rightarrow \mathcal{R}^+$. A weighted graph $(G_1,w_1)$ is called the weighted subgraph of $(G,w)$ if $G_1$ is a subgraph of $G$ and $w_1=w|_{E(G_1)}$. For simplicity of notations, we denote the weighted subgraph by $(G_1,w)$ instead of $(G_1,w_1)$.
\end{defi}

Note that unweighed graphs are just the special case where all the edge weights are equal to 1. Let $(G,w)$ be a weighted graph, and $\widetilde{G}$ the universal cover of $G$. By the covering map $\varphi\colon V(\widetilde{G})\rightarrow V(G);~\varphi((v_0,\dots,v_i))=v_i$,
the weighted function $w$ can lift in a natural way to a weighted function $\widetilde{w}\colon E(\widetilde{G})\rightarrow \mathcal{R}^+$, which is defined by $\widetilde{w}((v_0,\dots,v_{i-1})(v_0,\dots,v_{i-1},v_i))=w(v_{i-1}v_i)$. Thus we naturally get a weighted universal cover $(\widetilde{G},\widetilde{w})$ from
the weighted graph $(G,w)$. For simplicity of notation, we write $(\widetilde{G},w)$ instead of $(\widetilde{G},\widetilde{w})$.
\begin{defi}
For a weighted graph $(G,w)$ of order $n$. The adjacency matrix $A(G,w)$ of $(G,w)$ is defined by
$$(A(G,w))_{u,v}=\left\{
\begin{array}{ll}
  w(uv),& \text{if}~ uv\in E(G); \\
  0, & \text{otherwise}. \\
 \end{array}
 \right.$$
The weighted spectral radius of $G$ is the spectral radius of $A(G,w)$, and denoted by $\lambda_1(G,w)$. Order the eigenvalues of $A(G,w)$ as $\lambda_1(G,w)\ge \lambda_2(G,w)\ge \dots \ge\lambda_n(G,w)$.
\end{defi}

\section{Proof of Theorem \ref{thm1} and Corollaries}
The proof of Theorem \ref{thm1} uses an old idea of constructing a weighted test function via non-backtracking walks (see e.g. \cite{Jiang2019,Chung2016,Srivastava2018}), and also via the eigenvector of a path (see e.g. \cite{Jiang2019}).

~

\begin{Proof of 1}
For $e=(v_0,v_1)\in W_1$, let $T_e$ be the component of $\widetilde{G}(v_0,r+1)-(v_0)$ containing the vertex $e$, and it is a tree. Let $T$ be the disjoint union of a class of graphs $\{T_e\}_{e\in W_1}$, and thus $T$ is a forest. The vertex set of $T$ is $\bigcup_{i=1}^{r+1}W_i$, where $W_i$ is defined as the set of all non-backtracking walks of length $i$ on $G$.
By regarding every vertex $(v_0,v_1,\dots,v_i)$ in $T_e$, where $1\le i\le r+1$, as the vertex $(v_1,\dots,v_i)$ in $\widetilde{G}(v_1,r)$, we observe that $(T_e,w)$ is a weighted subgraph of $(\widetilde{G}(v_1,r),w)$.
By monotonicity of weighted spectral radius, $\lambda_1(\widetilde{G}(v_1,r),w)\ge\lambda_1(T_e,w)$ for  $e=(v_0,v_1)\in W_1$. Since $\lambda_1(T,w)=\max\{\lambda_1(T_e,w)\colon e\in W_1\}$, there exists a vertex $e^{\ast}$ such that $\lambda_1(T,w)=\lambda_1(T_{e^{\ast}},w)$. It derives that there exists a vertex $v_1^{\ast}$, the terminal vertex of $e^{\ast}$, such that $\lambda_1(\widetilde{G}(v_1^{\ast},r),w)\ge\lambda_1(T,w)$. It suffices to prove that for any function $g\colon V\rightarrow \mathcal{R}^+$,
$$\lambda_1(T,w)\ge 2\cos(\frac{\pi}{r+2})\cdot\frac{\sum_{(v_1,v_2)\in W_1}\sqrt{d(v_1)-1}~w(v_1v_2)~\sqrt{g(v_1)g(v_2)}}
{\sum_{v\in V(G)} g(v)d(v)}.$$

 Now we consider a Markov chain on $W_1$ as follows: the initial state $E_1$ is chosen from $W_1$ uniformly at random, and if the current state $E_i=(v_{i-1},v_i)$ is given, the next state $E_{i+1}$ will be chosen from $\{(v_i,v_{i+1})\in W_1\colon v_{i+1}\neq v_{i-1}\}$ uniformly at random. The transition matrix $P$ is
 $$
 P_{(u,v),(w,z)}=\left\{
                   \begin{array}{lll}
                     \frac{1}{d(v)-1}, & \text{if }v=w \text{ and } z\neq u;\hfill& \\
                     0,&  \text{otherwise}.\hfill&
                   \end{array}
                 \right.
 $$
 We can attach $E_1,\dots, E_i$ one by one to form a non-backtracking walk on $G$ of length $i$, which is denoted by the random variables $Y_i=(X_0,X_1\dots,X_i)$.

It is known that $\lambda=2\cos(\frac{\pi}{r+2})$ is the spectral radius of the path $P_{r+1}$ on $r+1$ vertices. Let $(x_1,\dots,x_{r+1})\in \mathcal{R}^{r+1}$ be a positive eigenvector of $P_{r+1}$ associated with $\lambda$. By the Rayleigh principle, it follows that
\begin{equation}\label{eq1}
 \sum_{i=2}^{r+1}2x_{i-1}x_i=\lambda\cdot\sum_{i=1}^{r+1}x_i^{2}.
\end{equation}

Define the vector
\begin{equation}
f\colon \bigcup_{i=1}^{r+1}W_i\rightarrow \mathcal{R}; ~ f(\omega)=x_i\sqrt{g(v_i) {\rm Pr}(Y_i=\omega)}
\end{equation}
 for $\omega=(v_0,v_1,\dots,v_i)\in W_i$, where $g\colon  V(G)\rightarrow\mathcal{R}^+$ is a fixed vertex weight function. Let $A(w)$ be the adjacency matrix of the weighted forest $(T,w)$. For $\omega=(v_0,\dots,v_{i-1},v_i)$, let $\omega^-= (v_0,\dots,v_{i-1})$. By simple calculations, we have
\begin{align}
\langle f,f \rangle&=\sum_{i=1}^{r+1}\sum_{\omega\in W_i} f(\omega)^2=\sum_{i=1}^{r+1}x_i^2\sum_{\omega\in W_i} g(v_i) {\rm Pr}(Y_i=\omega),\label{eq2}\\
\langle f,A(w)f \rangle&=\sum_{i=2}^{r+1}\sum_{\omega\in W_i} 2f(\omega^-)f(\omega)\cdot w(\omega^-\omega)\nonumber\\
 &=\sum_{i=2}^{r+1}2x_{i-1}x_i\sum_{\omega\in W_i}w(\omega^-\omega)\sqrt{g(v_{i-1})g(v_i)}\sqrt{{\rm Pr}(Y_{i-1}=\omega^-)\cdot {\rm Pr}(Y_i=\omega)}.\label{eq3}
\end{align}
For simplicity of notation, let
\begin{align*}
  I_i & =\sum_{\omega\in W_i} g(v_i) {\rm Pr}(Y_i=\omega), \\
  J_i & =\sum_{\omega\in W_i}w(\omega^-\omega)\sqrt{g(v_{i-1})g(v_i)}\sqrt{{\rm Pr}(Y_{i-1}=\omega^-)\cdot {\rm Pr}(Y_i=\omega)}.
\end{align*}

In order to complete the proof, it suffices to simplify $I_i$ and $J_i$ for every $i$.
Firstly, we have
\begin{align}\label{eq4}
 I_i={\rm E}\left[ g(X_i)\right]=\sum_{v\in V(G)}g(v)\cdot{\rm Pr}(X_{i}=v).
\end{align}
for $i\ge 1$.
Secondly, by the Markov property, for $i\ge 2$ and $\omega=(v_0,v_1,\dots,v_i)\in W_i$,
\begin{align*}
  {\rm Pr}(Y_{i-1}=\omega^-)
  =\frac{{\rm Pr}(Y_i=\omega)}{{\rm Pr}(E_i=(v_{i-1},v_i)|E_{i-1}=(v_{i-2},v_{i-1}))}
  =(d(v_{i-1})-1){\rm Pr}(Y_i=\omega).
\end{align*}
Note that $w(\omega^-\omega)$ is defined to be $w(v_{i-1}v_i)$. Thus it follows that for $i\ge 2$,
\begin{align}\label{eq5}
J_i&=\sum_{\omega=(v_0,v_1,\dots,v_i)\in W_i}\sqrt{d(v_{i-1})-1}~w(v_{i-1}v_i)\sqrt{g(v_{i-1})g(v_i)}~{\rm Pr}(Y_{i}=\omega)\nonumber\\
&={\rm E}\left[\sqrt{d(X_{i-1})-1}~w(X_{i-1}X_i)\sqrt{g(X_{i-1})g(X_i)}\right]\nonumber\\
&=\sum_{(v_1,v_2)\in W_1}\sqrt{d(v_1)-1}~w(v_1v_2)~\sqrt{g(v_1)g(v_2)}~{\rm Pr}(X_{i-1}=v_1,X_{i-2}=v_2).
\end{align}

Now we focus on the probabilities in \eqref{eq4} and \eqref{eq5}. Since the minimum degree of $G$ is at least $2$, the Markov chain has no absorbing states. And it is easy to see that the uniform distribution $x=(\frac{1}{|W_1|},\dots,\frac{1}{|W_1|})$ on $W_1$ is a stationary distribution of the Markov chain, that is $x=xP^{i-1}$ for $i\ge 1$, where $P$ is the transition matrix.
Thus we have ${\rm Pr}(E_i=e)=1/|W_1|$ for $i\ge 1$ and $e\in W_1$, which derives that for $i\ge 1$ and $v\in V(G)$,
 $${\rm Pr}(X_i=v)=\sum_{e\sim v}{\rm Pr}(E_i=e)=\frac{d(v)}{|W_1|},$$
where $e\sim v$ denotes ranging over all the edges incident to $v$. Hence for $ i\ge 2$ and $(v_1,v_2)\in W_1$,
\begin{align*}
{\rm Pr}(X_{i-1}=v_1,X_i=v_2)&=\sum_{\{u\colon (u,v_1)\in W_1,u\neq v_2\}}{\rm Pr}(E_{i-1}=(u,v_1),E_i=(v_1,v_2))\\
&=\sum_{\{u\colon (u,v_1)\in W_1,u\neq v_2\}}{\rm Pr}(E_i=(v_1,v_2)|E_{i-1}=(u,v_1)){\rm Pr}(E_{i-1}=(u,v_1))\\
&=\sum_{\{u\colon (u,v_1)\in W_1,u\neq v_2\}}\frac{1}{(d(v_1)-1)|W_1|}=\frac{1}{|W_1|}.
\end{align*}
By rewriting ${\rm Pr}(X_i=v)$ in \eqref{eq4} and ${\rm Pr}(X_{i-1}=v_1,X_i=v_2)$ in \eqref{eq5}, we can simplify $I_i$ and $J_i$, and substitute them to \eqref{eq2} and \eqref{eq3} to obtain
\begin{align*}
\langle f,f \rangle&=\sum_{i=1}^{r+1}x_i^2\sum_{v\in V(G)} g(v)\frac{d(v)}{|W_1|},\\
 \langle f,A(w)f \rangle&=\sum_{i=2}^{r+1}2x_{i-1}x_i\sum_{(v_1,v_2)\in W_1}\sqrt{d(v_1)-1}~w(v_1v_2)~\sqrt{g(v_1)g(v_2)}\frac{1}{|W_1|}.
\end{align*}
Finally, combining \eqref{eq1}, the equalities above, and the Rayleigh principle, we obtain
$$\lambda_1(T,w)\ge\frac{\langle f,A(w)f \rangle}{\langle f,f \rangle}= 2\cos(\frac{\pi}{r+2})\frac{\sum_{(v_1,v_2)\in W_1}\sqrt{d(v_1)-1}~w(v_1v_2)~\sqrt{g(v_1)g(v_2)}}
{\sum_{v\in V(G)} g(v)d(v)},$$
for any vertex weight function $g$, and complete the proof.

\end{Proof of 1}

If the weight function $w_0$ is defined by $w_0(uv)=(d(u)d(v))^{-\frac 12}$ for every edge $uv$ of a graph $G$, then let $g(v)=d(v)$ for every vertex $v$ of $G$ in Theorem \ref{thm1}, and we derive the following corollary.

\begin{coro}\label{coro1}
If a connected weighted graph $(G,w_0)$ with edge weight $w_0(uv)=(d(u)d(v))^{-\frac 12}$ has minimum degree at least $2$, then for any natural number $r$ with $r\ge 1$, there exists a vertex $v$ of $G$ such that
$$\lambda_1(\widetilde{G}(v,r),w_0)\ge \frac{2\sum_{u\in V(G)}d(u)\sqrt{d(u)-1}}{\sum_{u\in V(G)} d(u)^2}\cos\left(\frac{\pi}{r+2}\right).$$
\end{coro}
 Since the weighted unraveled ball $(\widetilde{G}(v,r),w_0)$ is a weighted induced subgraph of the weighted universal cover $(\widetilde{G},w_0)$, it follows that
 $\lambda_1(\widetilde{G},w_0)\ge\lambda_1(\widetilde{G}(v,r),w_0)$ by the monotonicity of weighted spectral radius.
Thus we can obtain a lower bound on $\lambda_1(\widetilde{G},w_0)$ by letting $r$ go to infinity in Corollary \ref{coro1}.

\begin{coro}\label{coro11}
If a connected weighted graph $(G,w_0)$ with edge weight $w_0(uv)=(d(u)d(v))^{-\frac 12}$ has minimum degree at least $2$, then the weighted spectral radius of its universal cover satisfies
$$\lambda_1(\widetilde{G},w_0)\ge \frac{2\sum_{u\in V(G)}d(u)\sqrt{d(u)-1}}{\sum_{u\in V(G)} d(u)^2}.$$
\end{coro}

\section{Weighted spectral radius of a ball}
For a weighted graph, {\it the weight of a closed walk} is the product of weights of all edges on the closed walk. The following result is well-known.

\begin{lem}\cite{Mohar1989}\label{lem2}
For any connected weighted graph $(G,w)$ (possibly infinite) and every vertex $v$ of $G$, the weighted spectral radius of $(G,w)$ is $$\lambda_1(G,w)=\limsup_{k\rightarrow\infty}\sqrt[2k]{t_{2k}^{(w)}(v)},$$
where $t_{2k}^{(w)}(v)$ is the total weight of all closed walks of length $2k$ from $v$ to itself in $G$.
\end{lem}
The following lemma establishes connections between the weighted spectral radius of a ball and its corresponding unraveled ball, and is an extension of \cite[Theorem 2.2]{Mohar2010}.

\begin{lem}\label{lem1}
For every vertex $v$ of a weighted graph $(G,w)$ and any natural number $r$ with $r\ge 1$,
$$\lambda_1(G(v,r),w)\ge \lambda_1(\widetilde{G}(v,r),w).$$
\end{lem}

\begin{Proof}
Recall that the vertex set of $\widetilde{G}(v,r)$ consists of all non-backtracking walks of length at most $ r$ starting at $v$.
By the covering map $\varphi\colon V(\widetilde{G})\rightarrow V(G);~\varphi((v_0,\dots,v_i))=v_i$, we can naturally construct a map $\sigma$, mapping a closed walk $(\omega=\omega_0,\dots,\omega_{2k}=\omega)$ of length $2k$ in $\widetilde{G}(v,r)$ to a closed walk $(v=v_0,\dots,v_{2k}=v)$ of length $2k$ in $G(v,r)$ for $k\ge 0$, where $v_j$ is the terminal vertex of $\omega_j$ for every $j$.

It is obvious that $\sigma$ is an injective map. In fact, since the covering map $\varphi$ is a local isomorphism, there exists an inverse map $\tau$ such that $\tau\sigma=\text{id}$. 
In addition, by naturally lifting the weight function $w$ to a weight function of $\widetilde{G}$, and the weights of walks are invariant under the map $\sigma$. Hence, the sum of the weights of closed walks of length $2k$ on $\widetilde{G}(v,r)$ is no more than the sum of the weights of closed walks of length $2k$ on $G(v,r)$. Therefore, Lemma \ref{lem2} gives that $\lambda_1(G(v,r),w)\ge \lambda_1(\widetilde{G}(v,r),w)$.

\end{Proof}

Combining Theorem \ref{thm1} and Lemma \ref{lem1}, we prove a lower bound on the weighted spectral radius of a ball of a graph with edge weight $w_0(uv)=(d(u)d(v))^{-\frac 12}$, which also has its own interest.

\begin{thm}\label{thm3}
Let $(G,w_0)$ be a weighted graph with edge weight $w_0(uv)=(d(u)d(v))^{-\frac 12}$. If the graph $G$ has average degree $d$ with $d\ge 2$ and second order average degree $\widetilde{d}$,
then for any natural number $r$ with $r\ge 1$, there exists a vertex $v$ of $G$ such that
\begin{equation}
 \lambda_1(G(v,r),w_0)
\ge \frac{2\sqrt{d-1}}{\widetilde{d}} \cos\left(\frac{\pi}{r+2}\right).
\end{equation}
\end{thm}
\begin{Proof}
Since $G$ may have vertices of degree 1, we cannot use Corollary \ref{coro1} directly. Instead, we consider the 2-core $H$ in $G$, the largest induced subgraph of $G$ with minimum degree at least $2$.
Observe that the 2-core $H$ can be obtained from $G$ by deleting vertices of degree 1 sequentially.
Since removing vertices of degree 1 from a graph of average degree at least cannot decrease its average degree, the 2-core $H$ is non-empty.
We decompose the proof into two parts.

\textbf{Case 1:} The 2-core $H$ is connected. Similarly with Corollary \ref{coro1}, by setting $w_0(uv)=(d_{G}(u)d_{G}(v))^{-\frac 12}$ and $g(v)=d_{G}(v)$ in Theorem \ref{thm1} we can derive that there exists a vertex $v$ of $H$ such that
$$\lambda_1(\widetilde{H}(v,r),w_0)
\ge \frac{2\sum_{u\in V(H)}d_{H}(u)\sqrt{d_{H}(u)-1}}{\sum_{u\in V(H)}d_{H}(u)d_G(u)}
\cos\left(\frac{\pi}{r+2}\right).$$
Since $(H(v,r),w_0)$ is a weighted subgraph of $(G(v,r),w_0)$, the monotonicity of weighted spectral radius and Lemma \ref{lem1} derive that
$$\lambda_1(G(v,r),w_0)
\ge\lambda_1(H(v,r),w_0)\ge\lambda_1(\widetilde{H}(v,r),w_0)
.$$
Combing the inequalities above, we only need to prove that
\begin{equation}\label{eq6}
 \frac{\sum_{u\in V(H)}d_{H}(u)\sqrt{d_{H}(u)-1}}{\sum_{u\in V(H)}d_{H}(u)d_G(u)}
 \ge \frac{\sqrt{d-1}}{\widetilde{d}}.
\end{equation}

Recall that the average degree of $H$ is at least $d$.
Note that $h(x)=x\sqrt{x-1}$ is a convex function for $x\ge 2$. It follows from Jensen's inequality that
\begin{equation}\label{eq6.5}
 \sum_{u\in V(H)}d_{H}(u)\sqrt{d_{H}(u)-1}
\ge\sum_{u\in V(H)}d_{H}(u)\sqrt{d-1}.
\end{equation}
Then it suffices to prove
\begin{equation}\label{eq7}
 \frac{\sum_{u\in V(H)}d_{H}(u)d_G(u)}
{\sum_{u\in V(H)}d_{H}(u)}
\le\widetilde{d}=\frac{\sum_{u\in V(G)}d_G(u)^2}{\sum_{u\in V(G)}d_G(u)},
\end{equation}
since \eqref{eq6.5} and \eqref{eq7} imply \eqref{eq6}.
Let $H^{\prime}$ be the spanning subgraph of $G$ composed of $H$ plus all isolated vertices in $V\setminus V(H)$. It is obvious that
\begin{equation}\label{eq8}
 \frac{\sum_{u\in V(H)}d_{H}(u)d_G(u)}
{\sum_{u\in V(H)}d_{H}(u)}
=\frac{\sum_{u\in V(G)}d_{H^{\prime}}(u)d_G(u)}
{\sum_{u\in V(G)}d_{H^{\prime}}(u)}.
\end{equation}

Recall that the 2-core $H$ can be obtained from $G$ by deleting vertices of degree 1 sequentially.
Thus we can recover $G$ by $H^{\prime}$, by sequentially adding some edges of $G$ in the opposite order, and each of the edges joins a non-isolated vertex and an isolated vertex in the current state.
Assume that there are $m$ edges to be added. For simplicity of notation, let $G^{(0)}=H^{\prime}$ and $G^{(m)}=G$.
In the $i$th step, assume that some edge $v_1v_2$ of $G$ is added to $G^{(i-1)}$, and $G^{(i)}$ is the resulting graph. 
Note that $\frac xy\le\frac{x+1}{y+2}$ for $y\ge 2x$. It follows that for all $i$,
$$\frac{d_{G^{(i)}}(v_j)}{\sum_{u\in V(G)}d_{G^{(i)}}(u)}\le\frac{d_{G^{(i)}}(v_j)+1}{\sum_{u\in V(G)}d_{G^{(i)}}(u)+2}
=\frac{d_{G^{(i+1)}}(v_j)}{\sum_{u\in V(G)}d_{G^{(i+1)}}(u)},~\forall ~j=1,2,$$
$$\frac{d_{G^{(i)}}(z)}{\sum_{u\in V(G)}d_{G^{(i)}}(u)}\le\frac{d_{G^{(i)}}(z)}{\sum_{u\in V(G)}d_{G^{(i)}}(u)+2}
=\frac{d_{G^{(i+1)}}(z)}{\sum_{u\in V(G)}d_{G^{(i+1)}}(u)},~\forall ~z\in V(G)\setminus\{v_1,v_2\}.$$
Consequently, we have
\begin{equation}\label{eq9}
 \frac{d_{H^{\prime}}(z)}{\sum_{u\in V(G)}d_{H^{\prime}}(u)}\le
\frac{d_{G^{(1)}}(z)}{\sum_{u\in V(G)}d_{G^{(1)}}(u)}\le\dots\le
\frac{d_{G}(z)}{\sum_{u\in V(G)}d_{G}(u)},~\forall~ z\in V(G).
\end{equation}
Therefore, \eqref{eq7} follows by \eqref{eq8} and \eqref{eq9}.

\textbf{Case 2:}  The 2-core $H$ is disconnected. Now $G$ is also disconnected. Assume that the disconnected graph $G$ is composed of $G_1,\dots,G_t$.
Let $H_i$ be the connected 2-core of $G_i$ with vertex set $V_i$ for $1\le i\le t$.
By the same argument with Case 1, there exists a vertex $u_i$ of $H_i$ such that
$$\lambda_1(G(u_i,r),w_0)
\ge \frac{2\sum_{u\in V_i}d_{H}(u)\sqrt{d_{H}(u)-1}}{\sum_{u\in V_i}d_{H}(u)d_G(u)}
 \cos\left(\frac{\pi}{r+2}\right)=\colon \frac{M_i}{N_i} \cos\left(\frac{\pi}{r+2}\right).$$
One can simply verify that there exists an $i_0\in[1,t]$ such that
$$\frac{M_1+\dots+M_t}{N_1+\dots+N_t}\le \max\left\{\frac{M_1}{N_1},\dots,\frac{M_t}{N_t}\right\}=\frac{M_{i_0}}{N_{i_0}}.$$
Thus there exists a vertex $u_{i_0}$ of $H_{i_0}$ such that
$$\lambda_1(G(u_{i_0},r),w_0)
\ge \frac{2\sum_{u\in V(H)}d_{H}(u)\sqrt{d_{H}(u)-1}}{\sum_{u\in V(H)}d_{H}(u)d_G(u)}
 \cos\left(\frac{\pi}{r+2}\right).$$
Using the same argument with Case 1, we can prove that
there exists a vertex $u_{i_0}$ of $G$ such that
$$
\lambda_1(G(u_{i_0},r),w_0)\ge \frac{2\sqrt{d-1}}{\widetilde{d}} \cos\left(\frac{\pi}{r+2}\right).
$$
\end{Proof}

\section{Proof of Theorem \ref{thm2}}
Recall that a graph $G$ is $(r,d,\widetilde{d},s)$-robust if for the induced subgraph of $G$ obtained by sequentially deleting any $s$ ball of radius $r$,
its average degree is at least $d$, and its second order average degree is at most $ \widetilde{d}$.
Before proving Theorem \ref{thm2}, we provide a lower bound on the
$s$-th largest weighted eigenvalue of a graph.
\begin{lem}\label{lem4}
Let $(G,w_0)$ be a weighted graph with edge weight $w_0(uv)=(d(u)d(v))^{-\frac 12}$. Given real numbers, $\widetilde{d}\ge d\ge 2$, and natural numbers, $r\ge 1$ and $s\ge 2$, if $G$ is $(r,d,\widetilde{d},s-1)$-robust,
then
$$\lambda_s(G,w_0)\ge \frac{2\sqrt{d-1}}{\widetilde{d}}\cos\left(\frac{\pi}{r+1}\right).$$
\end{lem}
\begin{Proof}
We show how to sequentially construct $G_1,\dots,G_s$, a collection of pairwise disjoint weighted induced subgraphs of $G$ such that no vertex in $V(G_i)$ is adjacent to a vertex of $V(G_j)$ for $i\neq j$, and
$\lambda_1(G_{i},w_0)\ge\frac{2\sqrt{d-1}}{\widetilde{d}} \cos\left(\frac{\pi}{r+1}\right)$ holds for $1\le i\le s$.

For simplicity of notation, let $G^{(0)}=G$. For $1\le i\le s,$
assume that we have obtained $ G^{(i-1)}$, an induced subgraph of $G$, by sequentially deleting $i-1$ balls of radius $r$ during the previous steps. In the $i$th step, let $H_i$ be the graph obtained by sequentially deleting any $s-i$ balls of radius $r$ from $G^{(i-1)}$. Since $G$ is $(r,d,\widetilde{d},s-1)$-robust,
the average degree of $H_i$ is at least $d$, and the second order average degree of $H_i$ is at most $\widetilde{d}$.
By Theorem \ref{thm3}, there exists a vertex $v_i$ of $H_i$ such that
$$\lambda_1(H_i(v_i,r-1),w_0)
\ge\frac{2\sqrt{d-1}}{\widetilde{d}} \cos\left(\frac{\pi}{r+1}\right).$$
Since $H_i$ is a subgraph of $G^{(i-1)}$, $H_i(v_i,r-1)$ is a subgraph of $G^{(i-1)}(v_i,r-1)=\colon G_{i}$.
Thus as desired we obtain
$$\lambda_1(G_{i},w_0)\ge\lambda_1(H_i(v_{i},r-1),w_0)\ge\frac{2\sqrt{d-1}}{\widetilde{d}} \cos\left(\frac{\pi}{r+1}\right).$$
Let $G^{(i)}$ be the induced subgraph of $G^{(i-1)}$ obtained by deleting a ball $G^{(i-1)}(v_i,r)$. Turn to the next step, until we get $G_1,\dots,G_s$.

For $1\le i\le s$, let $A_i(w_0)$ be the adjacency matrix of $(G_i,w_0)$. Let $A(w_0)$ be the adjacency matrix of $(G,w_0)$.
Note that $A_i(w_0)$ is just
the principal sub-matrix of $A(w_0)$ corresponding to $G_i$, but not equal to $I-\mathcal{L}(G_i)$.
Additionally, let $f_i$ be a positive unit eigenvector of $A_i(w_0)$ associated with $\lambda_1(G_i,w_0)$.
We can define a vector $g_i\colon V(G)\rightarrow\mathcal{R}$ by
$$g_i(u)=
\left\{
      \begin{array}{cc}
           f_i(u),\hfill & \text{if } u\in V(G_i); \hfill\\
           0,\hfill&\text{otherwise},\hfill\\
         \end{array}
       \right.
$$
for every $i$ with $1\le i\le s$. It is obvious that $V(G_i)\cap V(G_j)=\emptyset$ for $ i\neq j$, so $\{g_1,\dots,g_s\}$
is a set of orthonormal vectors. Define $W_0$ as an $s$-dimension vector space spanned by $\{g_1,\dots,g_s\}$.
By the Rayleigh Principle, it follows that
$$\lambda_s(G,w_0)=\max_{\text{dim}W=s}\min_{f\in W}\frac{\langle f,A(w_0) f\rangle}{\langle f, f\rangle}\ge\min_{f\in W_0}\frac{\langle f,A(w_0) f\rangle}{\langle f, f\rangle}.$$
It is obvious that $uv\notin E$ for $u\in V(G_i)$ and $v\in V(G_j)$ with $i\neq j$. Taking any non-zero element $f=c_1g_1+\dots+c_sg_s$ in $W_0$, we have
\begin{align*}
 \frac{\langle f,A(w_0) f\rangle}{\langle f, f\rangle}
 &=\frac{c_1^2\langle f_1,A_1(w_0) f_1\rangle+\dots+c_s^2\langle f_s,A_s(w_0) f_s\rangle}{c_1^2\langle f_1, f_1\rangle+\dots+c_s^2\langle f_s, f_s\rangle} \\
 & =\frac{c_1^2\lambda_1(G_1,w_0)\langle f_1,f_1\rangle+\dots+c_s^2\lambda_1(G_s,w_0)\langle f_s, f_s\rangle}{c_1^2\langle f_1, f_1\rangle+\dots+c_s^2\langle f_s, f_s\rangle}\\
 &\ge\frac{2\sqrt{d-1}}{\widetilde{d}} \cos\left(\frac{\pi}{r+1}\right).
\end{align*}
Finally, we complete the proof by
$$\lambda_s(G,w_0)\ge\min_{f\in W_0}\frac{\langle f,A(w_0) f\rangle}{\langle f, f\rangle}\ge\frac{2\sqrt{d-1}}{\widetilde{d}} \cos\left(\frac{\pi}{r+1}\right).$$
\end{Proof}

\begin{Proof of 2}
Recall that the normalized Laplacian matrix $\mathcal {L}$ of $G$ is defined to be $I-D^{-\frac 12}AD^{-\frac 12}$, where $D$ is the diagonal degree matrix of $G$ and $A$ is the adjacency matrix of $G$.
If we consider a weighted graph $(G,w_0)$ with edge weight $w_0(uv)=(d(u) d(v))^{-\frac 12}$, then we have $\mathcal {L}=I-A(G,w_0)$.
By applying Lemma \ref{lem4}, we can obtain the desired upper bound on $\mu_s(G)=1-\lambda_s(G,w_0)$.

\end{Proof of 2}

\subsection*{Acknowledgements}
The authors would be grateful to the referees for their valuable suggestions and comments which make a great improvement of the manuscript.


\begin{thebibliography}{99}

\bibitem{alon1986}
N. Alon,
Eigenvalues and expanders,
Theory of computing (Singer Island, Fla., 1984). Combinatorica 6 (1986), no. 2, 83--96.


\bibitem{angel2015}
O. Angel, J. Friedman, S. Hoory,
The non-backtracking spectrum of the universal cover of a graph,
Trans. Amer. Math. Soc. 367 (2015), 4287--4318.


\bibitem{Bordenave2019}
C. Bordenave, S. Coste,
Graphs with prescribed local neighborhoods of their universal coverings,
J. Combin. Theory Ser. B 138 (2019), 196--205.

\bibitem{chung1997} F. R. K. Chung,  Spectral graph theory, CBMS Regional Conference Series in Mathematics, 92. Published for the Conference Board of the Mathematical Sciences, Washington, DC; by the American Mathematical Society, Providence, RI, 1997.


\bibitem{Chung2016}
F. Chung,
A generalized Alon--Boppana bound and weak Ramanujan graphs,
Electron. J. Combin. 23 (2016), no. 3, Paper 3.4, 20 pp.

\bibitem{Hoory2005}
S. Hoory,
A lower bound on the spectral radius of the universal cover of a graph,
J. Combin. Theory Ser. B 93 (2005), no. 1, 33--43.

\bibitem{Jiang2019}
Z. Jiang,
On spectral radii of unraveled balls,
J. Combin. Theory Ser. B 136 (2019), 72--80.

\bibitem{Leighton1982}
F. T. Leighton,
Finite common coverings of graphs,
 J. Combin. Theory Ser. B 33 (1982), no. 3, 231--238.

\bibitem{Mohar2010}
B. Mohar,
A strengthening and a multipartite generalization of the Alon--Boppana--Serre theorem,
Proc. Amer. Math. Soc. 138 (2010), no. 11, 3899--3909.

\bibitem{Mohar1989}
B. Mohar, W. Woess,
A survey on spectra of infinite graphs,
Bull. London Math. Soc. 21 (1989), no. 3, 209--234.

\bibitem{Nilli1991}
A. Nilli,
On the second eigenvalue of a graph,
Discrete Math. 91 (1991), no. 2, 207--210.




\bibitem{Srivastava2018}
N. Srivastava, L. Trevisan,
An Alon--Boppana type bound for weighted graphs and lower bounds for spectral sparsification,
Proceedings of the Twenty-Ninth Annual ACM-SIAM Symposium on Discrete Algorithms, 1306--1315, SIAM, Philadelphia, PA, 2018.

\bibitem{Young2022}
S. J. Young,
The weighted spectrum of the universal cover and an Alon--Boppana result for the normalized Laplacian,
J. Comb. 13 (2022), no. 1, 23--40.
\end{thebibliography}
\end{document}